\documentclass[a4paper,11pt]{amsart}

\usepackage{amsthm, amsmath, latexsym, amsfonts, graphics, url}
\usepackage[all]{xy}
\usepackage{amssymb}

\title[Guido Castelnuovo]{Guido Castelnuovo and his heritage: \\
geometry, combinatorics, teaching}
\author{Claudio Fontanari}

\email{claudio.fontanari@unitn.it}\curraddr{
{\sc Dipartimento di Matematica \\  Universit\`a degli Studi di Trento\\
Via Sommarive 14 \\ 38123 Trento \\ Italy.}}

\thanks{
This research project was partially supported by GNSAGA of INdAM and by PRIN 2017 ``Moduli Theory and Birational Classification".\\
}

\begin{document}

\begin{abstract}
We approach Guido Castelnuovo's intellectual world by focusing on a trilogy of papers published 
in 1889 and by drawing a few remarks about Castelnuovo's scientific interests and attitudes. 
\end{abstract}

\maketitle

\section{Introduction}

For an account of the life and works of Guido Castelnuovo (1965--1952) we refer to \cite{A}, which also provides 
a fascinating insight into his familiar milieu, and to the biographical sketch \cite{G}. Here instead we are going to 
focus our attention on three early papers, published by Castelnuovo in 1889, which are specifically devoted to 
enumerative geometry and algebraic curves, but shed light more generally on his intellectual world. 

Indeed, Castelnuovo's contribution to the development of Brill-Noether theory is well-known, 
a terse presentation in modern terms is provided for instance in \cite{HM}, p. 242:

\begin{quote}

Let's now examine this history a bit more closely. To begin with, Brill and Noether asserted the truth 
of the theorem based, apparently, on a naive dimension count bolstered by the calculation of examples
in low genus. 

Exactly how the desired variational element might enter into the proofs was first suggested 
by Castelnuovo. His goal was not to establish any of the present theorems. Rather, he assumed the
statement of the Brill-Noether theorem and applied it to compute the {\it number} of $g^r_d$'s on a 
general curve in the case $\rho = 0$, when we expect it to be finite. 

To do this, Castelnuovo looked not at any smooth curve of genus $g$, but at a $g$-nodal curve $C_0$: 
that is, a rational curve with $g$ nodes $r_1, \ldots , r_g$ obtained by identifying $g$ pairs of points 
$(p_i, q_i)$ on $\mathbb{P}^1$. 

Any $g^r_d$ on $C_0$, Castelnuovo reasoned, would pull back to a $g^r_d$ on $\mathbb{P}^1$, which could 
then be represented as the linear series cut out on a rational normal curve $C \cong \mathbb{P}^1 
\hookrightarrow \mathbb{P}^d$ of degree $d$ by those hyperplanes containing a fixed $(d-r-1)$-plane 
$\Lambda \subset \mathbb{P}^d$. 

The condition that the $g^r_d$ on $\mathbb{P}^1$ be the pullback of one on $C_0$ was simply that 
every divisor of the $g^r_d$ containing $p_i$ should contain $q_i$ as well; in other words,
$\Lambda$ should meet each of the chords $\overline{p_i q_i}$ to $C$ in $\mathbb{P}^d$. 

The number of $g^r_d$'s on a general curve of genus $g$ was thus, according to Castelnuovo, the
number of $(d-r-1)$-planes in $\mathbb{P}^d$ meeting each of $g$ lines, a problem
in Schubert calculus that Castelnuovo went on to solve (to obtain the
correct value for the number of $g^r_d$'s on a general curve).

It was Severi who first pointed out, some twenty years later, that Castelnuovo's computation 
might serve as the basis of a proof of the Brill-Noether statement.

\end{quote}

Our plan is to give a closer look at the original papers, by providing in Section \ref{first} 
an English translation of some key passages, and by drawing in Section \ref{second} 
a series of remarks based on Castelnuovo's own statements. 

This material has been presented at the conference TiME 2019 (Levico Terme, September 2--6, 2019). 
The author is grateful to the organizers for their kind invitation and to Enrico Arbarello for enlightening 
conversations about Guido Castelnuovo and the genuine meaning of his heritage. 

\section{A trilogy of papers} \label{first}

G. Castelnuovo, {\it Una applicazione della geometria enumerativa 
alle curve algebriche}, Rendiconti del Circolo Matematico di Palermo, t. III, 1889
(\cite{C}, III, pp. 45--53).

\begin{quote}
In questo lavoro ci proponiamo due fini: esporre un metodo utile in molte ricerche della teoria delle curve; 
presentare alcune formole che ci sembrano notevoli e in se stesse, e per le loro conseguenze. A queste formole 
noi siamo giunti applicando il principio della conservazione del numero a curve degeneri.
\end{quote}

\vspace{0.1cm}

\begin{quote}
In this paper we have a twofold goal: to discuss a method useful in several researches in curve theory; 
to present some formulas we believe remarkable both in themselves and for their consequences. 
We came to these formulas by applying the principle of conservation of the number to degenerate curves. 
\end{quote}

\vspace{0.3cm}

G. Castelnuovo, {\it Numero degli spazi che segano pi\`u rette in uno spazio 
ad $n$ dimensioni}, Rendiconti della R. Accademia dei Lincei, s. IV, vol. V, 4 agosto 1889
(\cite{C}, IV, pp. 55--64).

\begin{quote}
Fra le questioni che appartengono alla Geometria Enumerativa, va notata per la sua importanza algebrica 
e geometrica la seguente: Quanti sono gli spazi ad $s$ dimensioni che soddisfanno a pi\`u condizioni fondamentali 
date in uno spazio ad $n$ dimensioni? Naturalmente le condizioni si suppongono tali da rendere 
determinato il problema.

(...) 

Un caso particolare interessante \`e il seguente: 

{\it Date $hk$ rette in uno spazio a $\{(k+1)(h-1)\}$ dimensioni, il numero degli spazi a $\{k(h-1)-1 \}$
dimensioni che segano in punti queste rette \`e uguale a
$$
\frac{1!2!3! \ldots (h-1)! 1!2!3! \ldots (k-1)!(hk)!}{1!2!3! \ldots (h+k-1)!}.
$$}
\end{quote}

\vspace{0.1cm}

\begin{quote}
Among the questions pertaining to Enumerative Geometry, we point out the following one, due to its both algebraic and geometric 
relevance: How many are the $s$-dimensional spaces satisfying several fundamental conditions given in an $n$-dimensional space?
Obviously, we suppose these conditions are such that the problem is made determinate. 

(...)

A particular interesting case is the following: 

{\it Given $hk$ lines in a $\{(k+1)(h-1)\}$-dimensional space, the number of $\{k(h-1)-1 \}$-dimensional spaces
cutting these lines in points is equal to
$$
\frac{1!2!3! \ldots (h-1)! 1!2!3! \ldots (k-1)!(hk)!}{1!2!3! \ldots (h+k-1)!}.
$$}

\end{quote}

\vspace{0.3cm}

G. Castelnuovo, {\it Numero delle involuzioni razionali giacenti sopra una curva di dato 
genere}, Rendiconti R. Accademia dei Lincei, s. IV, vol. V, 1 settembre 1889
(\cite{C}, V, pp. 65--68).

\begin{quote}
\`E noto che sopra una curva di genere $p$ con moduli generali esistono delle {\it serie lineari} 
$g^{(q)}_m$ (involuzioni razionali di $\infty^q$ gruppi di $m$ punti) in numero finito, quando sia 
\begin{equation}\label{one}
m-q=(p-m+q)q;
\end{equation} 
{\it quante sono queste $g^{(q)}_m$?} (...) Noi ci proponiamo di risolvere il problema in tutta la sua generalit\`a 
approfittando del seguente concetto di geometria enumerativa, che ci serv\`i in altra 
occasione\footnote{V. la memoria III, \emph{Una applicazione della geometria enumerativa}. 
Dobbiamo riconoscere che nello stabilire questo concetto ci fondiamo pi\`u sulla intuizione (e su varie 
verificazioni), che sopra un vero ragionamento matematico. Alla dimostrazione si potr\`a forse arrivare 
considerando la curva in uno spazio superiore come intersezione parziale di pi\`u variet\`a e trattando 
algebricamente il problema degli spazi secanti; si troverebbe che il numero delle soluzioni \`e indipendente 
dalla posizione particolare delle variet\`a. Ma un ragionamento di tal natura non potr\`a farsi che quando 
la teoria delle curve negli spazi superiori sar\`a pi\`u completa. Ci permettiamo per\`o di approfittare 
di un principio non ancora dimostrato per risolvere un difficile problema, perch\'e crediamo che anche 
con simili tentativi si possa giovare alla scienza, quando si dichiari esplicitamente ci\`o che si ammette 
e ci\`o che si dimostra.}: 

"il numero (supposto finito) degli $[r]$ (spazi ad $r$ dimensioni) che segano in $\sigma$ punti una curva 
$C^n_p$ (di ordine $n$ e genere $p$) appartenente ad un $[s]$ non muta, o diventa infinito, quando 
alla curva data si sostituisca l'insieme di pi\`u curve, purch\'e l'ordine ed il genere della curva composta siano ancora risp. $n$ e $p$".

Noi useremo soltanto curve costituite da una curva semplice con pi\`u corde, e precisamente in luogo di una curva di genere $p$, 
considereremo una curva razionale insieme a $p$ delle sue corde.

(...)

La $C^{m+p}_p$ si scinda in una curva razionale $C^m_0$ appartenente ad $[m]$ ed in $p$ delle sue corde 
scelte ad arbitrio.

(...)

{\it Se \`e soddisfatta la (\ref{one}), il numero delle serie $g^{(q)}_m$ esistenti sopra una curva di genere $p$ uguaglia il numero degli spazi $[m-q-1]$ che segano $p$ rette di $[m]$.}

Ora quest'ultimo numero fu gi\`a determinato\footnote{V. la nostra nota IV: {\it Numero degli spazi che segano pi\`u rette in uno spazio ad $n$ dimensioni}, \S 10.};
se per semplicit\`a poniamo
$$
p-1-(m-q)=Q
$$
ossia: 
$$
m = p-1-(Q-q)
$$
e quindi per la (\ref{one})
$$
p = (q+1)(Q+1),
$$
si trova che il numero di cui si parla \`e dato da
$$
\frac{1!2!3! \ldots q! 1!2!3! \ldots Q!p!}{1!2!3! \ldots (q+Q+1)!}.
$$
\end{quote}

\vspace{0.1cm}

\begin{quote}
It is known that on a curve of genus $p$ with general moduli there exist a finite number of {\it linear series} 
$g^{(q)}_m$ (rational involutions of $\infty^q$ groups of $m$ points) whenever
\begin{equation}\label{onebis}
m-q=(p-m+q)q;
\end{equation} 
{\it how many are these $g^{(q)}_m$?} (...) Our goal is to solve this problem in complete generality
by exploiting the following notion of enumerative geometry, which we applied in another 
occasion\footnote{See the paper III, \emph{Una applicazione della geometria enumerativa}. 
We have to admit that in establishing this notion we rely more on intuition (an on several checks) 
than on a complete mathematical argument. Maybe the proof will be reached by considering the 
curve in a higher space as the partial intersection of several varieties and by treating algebraically 
the problem of secant spaces; we would find that the number of solutions in independent on the 
special position of the varieties. But such an argument could be carried out only when the theory 
of curves in higher spaces will be more developed. On the other hand, we are going to apply a 
still unproven principle in order to solve a difficult problem since we believe that such attempts
may be useful to the progress of science, provided one explicitly declare what is admitted and 
what is proven.}: 

"the number (assumed to be finite) of the $[r]$ ($r$-dimensional spaces) cutting $\sigma$ points 
on a curve $C^n_p$ (of degree $n$ and genus $p$) in a $[s]$ does not change, nor becomes 
infinite, when the given curve is replaced by the union of several curves, provided that the degree 
and the genus of the reducible curve are still $n$ e $p$, respectively".

We are going to use only curves made by an irreducible curve together with several chords, 
more precisely, instead of a curve of genus $p$ we will consider a rational curve together with 
$p$ of its chords. 

(...)

Let the $C^{m+p}_p$ be split into a rational curve $C^m_0$ in $[m]$ and into $p$ of its chords, arbitrarily chosen. 

(...)

{\it If (\ref{onebis}) is satisfied, the number of series $g^{(q)}_m$ laying on a curve of genus $p$ is equal to the number of space $[m-q-1]$ cutting $p$ lines of $[m]$.}

Now, this last number was already computed\footnote{See our note IV: {\it Numero degli spazi che segano pi\`u rette in uno spazio ad $n$ dimensioni}, \S 10.};
if for simplicity we set
$$
p-1-(m-q)=Q
$$
that is to say: 
$$
m = p-1-(Q-q)
$$
hence by (\ref{onebis})
$$
p = (q+1)(Q+1),
$$
we find that the number we are looking for is given by
$$
\frac{1!2!3! \ldots q! 1!2!3! \ldots Q!p!}{1!2!3! \ldots (q+Q+1)!}.
$$
\end{quote}

\vspace{0.3cm}

Guido Castelnuovo: {\it Memorie scelte}, Aggiunta alle Memorie III, IV e V, p. 69.

\begin{quote}

Nell'autunno 1888 C. Segre, che allora studiava le superficie rigate algebriche degli iperspazi, 
aveva ricondotto la determinazione del numero delle direttrici (curve unisecanti le generatrici) 
d'ordine minimo di una rigata al calcolo del numero degli spazi $S_{s-1}$ che contengono $s$ generatrici 
di una rigata di $S_s$, e mi propose il problema di stabilire questo numero. Di qua l'origine della 
Mem. III, che contiene al n. 4 il risultato richiesto, del quale ha tenuto conto il Segre 
("Mathematische Annalen", vol. 34, 1889). L'idea che mi ha permesso di raggiungere rapidamente 
questo e altri risultati consiste nel sostituire ad una curva irriducibile d'ordine $n$ e genere $p$ 
di un iperspazio, una curva composta di una curva d'ordine $n-1$ e di una retta unisecante o bisecante, 
secondo che quest'ultima curva ha genere $p$ o $p-1$. 

Questo \emph{principio di degenerazione}, che serve anche nella Mem. V, \`e 
semplicemente ammesso; la prima dimostrazione che lo spezzamento non altera i numeri richiesti fu data per via topologica 
(ricorrendo alle superficie di Riemann) da F. Klein in un suo corso del $2^o$ semestre 1892 (\emph{Riemannsche Fl\"achen}, 
II, pag. 110 e segg., lezioni litografate 1892). Per via algebrica occorre far vedere che la curva spezzata pu\`o esser riguardata 
come limite di una curva irriducibile variante entro un sistema continuo, ci\`o che, sotto ipotesi assai larghe, ha dimostrato F. Severi 
nelle \emph{Vorlesungen \"uber algebraische Geometrie}, Anhang G. (B. G. Teubner, Leipzig-Berlin, 1921); v. in particolare pag. 392.

\end{quote}

\vspace{0.1cm}

\begin{quote}

In autumn 1888 C. Segre, who was then investigating ruled algebraic surfaces in hyperspaces, 
had reduced the determination of the number of directrices (cubes unisecant the generatrices)
of minimal degree of a ruled surface to the computation of the number of spaces $S_{s-1}$
containing $s$ generatrices of a ruled surface in $S_s$, and he posed to me the problem of 
establishing this number. Hence the origin of paper III, containing at n. 4 the requested result, 
taken into account by Segre ("Mathematische Annalen", vol. 34, 1889). The idea that allowed me 
to quickly reach this and other results consists in substituting to an irreducible curve of degree $n$
and genus $p$ in a hyperspace, a reducible curve union of a curve of degree $n-1$ and of a line, 
either unisecant or bisecant according to the fact that this last curve has genus either $p$ or $p-1$. 

This \emph{principle of degeneration}, applied also in the paper V, is simply admitted; the first proof 
that the splitting does not change the required numbers was given through topology (by exploiting 
Riemann surfaces) by F. Klein in his course delivered on the second semester of 1892
(\emph{Riemannsche Fl\"achen}, II, p. 110 and ff., lithographed lectures 1892). To follow 
an algebraic way one needs to show that the reducible curve may be obtained as a limit 
of an irreducible curve varying in a continuous system, that is, under very mild assumptions, 
F. Severi proved in  \emph{Vorlesungen \"uber algebraische Geometrie}, Anhang G. 
(B. G. Teubner, Leipzig-Berlin, 1921); see in particular p. 392.

\end{quote}

\section{A few remarks} \label{second}

\subsection{Guido Castelnuovo and Corrado Segre}

We have seen that the starting point of Castelnuovo's research is an enumerative question posed to him 
by Corrado Segre. The role of the intense scientific dialogue between Segre and Castelnuovo in the 
flourishing of the Italian school of algebraic geometry should not be underestimated: in a letter to Amodeo 
sent on February 6, 1893, Castelnuovo mentions with longing the \emph{orge geometriche torinesi} 
(\emph{geometric orgies in Turin}). As pointed out by Enrico Arbarello (private communication):

\begin{quote}
When thinking of Castenuovo, and of his walks with Corrado Segre through the streets of Turin, 
during which they rapidly absorbed Riemann's point of view through the prism of their Italian taste, 
I can't help but think of how unfamiliar, abstract, and foreign it must have felt to them. I think that 
without being exposed to those revolutionary ideas, which freed him from the extrinsic world of Cremona
transformations, Castelnuovo might not have arrived at his contraction principle. Indeed, what I would celebrate 
in Castelnuovo is his absolute open-mindedness, his taste for adventure, and his eagerness to create and 
follow new paths, no matter how faint their trace.
\end{quote}

The mail correspondence between Segre and Castelnuovo witnesses their both scientific and personal 
deep fellowship. Here Segre comments Castel\-nuovo's results presented above:

\vspace{0.3cm}

C. Segre a G. Castelnuovo, 20 IX 1888 (\cite{AL})

\begin{quote}

\begin{flushright}
Torino, 20 IX 88
\end{flushright}

Carissimo Castelnuovo,

Alcuni dei teoremi che mi comunichi mi paiono veramente importanti. Importante l'idea di servirsi di curve 
di genere $p$ \underline{degeneri}. (...)

\end{quote}

\vspace{0.2cm}

C. Segre to G. Castelnuovo, 20 IX 1888

\begin{quote}

\begin{flushright}
Torino, 20 IX 88
\end{flushright}

My dearest Castelnuovo,

Some of the theorems you told me seem to me really important. Important the idea of exploiting 
\underline{degenerations} of curves of genus $p$. (...)

\end{quote}

\vspace{0.2cm}

Here instead Segre addresses touching words to Castelnuovo just after his move from Turin to Rome: 

\vspace{0.3cm}

C. Segre a G. Castelnuovo, 12 XI 1891 (\cite{AL})

\begin{quote}

\begin{flushright}
Torino, 12 XI 91
\end{flushright}

Mio carissimo, 

Ricevo la tua affettuosa lettera, e te ne ringrazio. Da Luned\`i tu mi manchi 
ed io sento vivamente questa lacuna. Tu accenni 
a quel po' di giovamento che hai potuto trarre in questi quattro anni dalla mia compagnia. 
Se ci\`o \`e vero, \`e pur vero che da te io ho avuto un completo ricambio, 
e che il tuo ingegno acuto, come la tua bont\`a di cuore 
m'han reso continuamente utili e piacevoli le tante ore che passavamo insieme (...) 
Tu m'hai fatto del bene, lo ripeto, non solo intellettualmente ma anche moralmente. 
Ed ora che tu mi manchi sento realmente un vuoto, che non sar\`a colmato da nessuno. (...)
Conservami sempre il tuo affetto. (...)

E ancora una volta un abbraccio affettuosissimo dal
\begin{flushright}
Tuo aff.mo C. Segre
\end{flushright}

\end{quote}

\vspace{0.1cm}

C. Segre to G. Castelnuovo, 12 XI 1891
(English version by Enrico Arbarello in \cite{A}, p. 21)

\begin{quote}

\begin{flushright}
Torino, 12 XI 91
\end{flushright}

My Dearest,

I am in receipt of your affectionate letter and I thank you for it.
Ever since Monday I miss you, and I feel this void very deeply. You
mention the bit of benefit that you might have been able to gain
from the past four years in my company. If that is indeed true, it is
also true that it was a completely even exchange, and that your
acute insight, as well as the goodness of your heart, have
continuously made the many hours I spent with you useful and
pleasant (...) 
You did me good, I repeat, not only intellectually,
but also morally. And now that you are missing, I really feel a void
which cannot be filled by anyone. (...) 
Keep me forever in your affection. (...)
And once again, a very big hug, from

\begin{flushright}
Your affectionate C. Segre
\end{flushright}

\end{quote}

\subsection{Geometry and probability}
We have seen that the number computed by Castelnuovo:
$$
\frac{1!2!3! \ldots q! 1!2!3! \ldots Q!p!}{1!2!3! \ldots (q+Q+1)!}
$$
presents a pronounced combinatorial flavour. It is well-known that Castel\-nuovo is the author 
of the first Italian treatise on the Calculus of Probability (1919) and this pioneering work \cite{CP} outside 
the realm of algebraic geometry is usually motivated by a gradual shift of his scientific interests towards 
applications of mathematics to natural and social phenomena (see for instance \cite{G}, p. 167). 
The above computation suggests a different (maybe complementary) explanation, by showing 
a much earlier Castelnuovo's taste for combinatorial manipulations analogous to the ones 
involved in discrete probability.  

\subsection{Guido ed Emma Castelnuovo}

We have seen that Castelnuovo claims to apply \emph{a still unproven principle in order to solve 
a difficult problem since... such attempts may be useful to the progress of science, 
provided one explicitly declare what is admitted and what is proven.} This intellectual habit, 
which is a very personal mixture of free open-mindedness and strict moral rigour, is typical 
of Castelnuovo's attitude to both mathematical research and teaching. His daughter, 
Emma Castelnuovo, consciously collects his heritage in her masterpiece volume \cite{EC}
devoted to school teaching of mathematics. In particular, two explicit quotations in her book 
are devoted to Guido Castelnuovo. The first one is on p. 5:

\vspace{0.1cm}

\begin{quote}
Guido Castelnuovo [espone] delle riflessioni fortemente indicative per un moderno insegnamento\footnote{G. Castelnuovo, 
{\it La scuola nei rapporti con la vita e la scienza moderna}, conferenza tenuta a Genova nel 1912 
in occasione del III Congresso della Mathesis, e riprodotta in {\it Archimede}, n. 2--3, 1962.}
e a cui, rilette a distanza di pi\`u di cinquanta anni, potrebbero ispirarsi oggi i compilatori dei 
programmi di matematica: 
{\it \`E questo il torto precipuo dello spirito dottrinario che invade
la nostra scuola. Noi vi insegnamo a diffidare dell'approssimazione, che \`e realt\`a, per 
adorare l'idolo di una perfezione che \`e illusoria.
Noi vi rappresentiamo l'universo come un edificio, le cui linee hanno una perfezione 
geometrica e ci sembrano sfigurate e annebbiate in causa del carattere grossolano dei 
nostri sensi, mentre dovremmo far comprendere che le forme incerte rivelateci dai sensi 
costituiscono la sola realt\`a accessibile, alla quale sostituiamo, per rispondere a certe 
esigenze del nostro spirito, una precisione ideale...}
\end{quote}

\vspace{0.1cm}

\begin{quote}
Guido Castelnuovo presents a few strongly suggestive remarks for a modern teaching\footnote{G. Castelnuovo, 
{\it La scuola nei rapporti con la vita e la scienza moderna}, talk held in Genova in 1912 
at the III Congress of Mathesis, and reproduced in {\it Archimede}, n. 2--3, 1962.}, 
from which school guidelines for mathematics could take inspiration today, more than fifty years later:
{\it This is the main fault of the dogmatic spirit invading our school. There we teach to distrust 
approximation, which is reality, in order to worship the idol of a perfection, which is illusion.  
There we represent the universe as a building, whose lines have a geometric perfection and look 
deformed and obfuscated because of the rough character of our senses, while we should make clear that 
the confused forms disclose to us by our senses are the only achievable reality, to which we replace, 
in order to satisfy certain needs of our spirit, an ideal precision...} 
\end{quote}

\vspace{0.1cm}

The second, even more amazing, quotation of Guido Castelnuovo appears on p. 157 of \cite{EC}:

\begin{quote}
Si dir\`a che \`e impossibile dare al bambino una nozione certa di funzione, che \`e pericoloso 
parlare del concetto di limite in termini vaghi, si dir\`a che quanto si insegna deve essere 
perfetto per non originare idee false che poi sarebbe difficile sradicare per sostituirle con 
appropriate definizioni. Mi torna alla mente quanto scriveva, nel lontano 1912, Guido 
Castelnuovo a questo proposito: {\it Ci\`o che si sa dal professore o dall'allievo -- mi fu detto --, 
sia pur limitato, ma deve sapersi perfettamente. Orbene, io sono uno spirito mite e tollerante; 
ma tutte le volte che questa frase mi fu obiettata, un maligno pensiero mi ha attraversato 
come un lampo la mente.
Oh, se potessi prendere in parola il mio interlocutore, e con un magico potere riuscissi a 
spegnere per un istante nel suo cervello tutte le cognizioni vaghe per lasciar sussistere 
soltanto ci\`o che egli sa perfettamente! Voi non immaginate mai quale miserando 
spettacolo potrei presentarvi! Ammesso pure che dopo una cos\`i crudele mutilazione 
qualche barlume rimanesse ancor nel suo intelletto, e di ci\`o fortemente dubito, 
somiglierebbe questo ad un gioco di fuochi folletti sperduti in tenebre profonde e 
sconfinate. La verit\`a \`e che noi nulla sappiamo perfettamente...}\footnote{G. Castelnuovo, 
{\it La scuola nei rapporti con la vita e la scienza moderna}, conferenza tenuta a Genova nel 1912 
in occasione del III Congresso della Mathesis, e riprodotta in {\it Archimede}, n. 2--3, 1962.} 
\end{quote}

\vspace{0.1cm}

\begin{quote}
One could object that it is impossible to provide a child with a rigorous notion of function, 
that it is dangerous to introduce the concept of limit in vague terms, one could object that 
what is taught has to be perfect in order to avoid misconceptions hard to eradicate and replace 
with appropriate definitions. It comes back to my mind what Guido Castelnuovo wrote about this, 
long ago in 1912: {\it What is known by the teacher or by the pupil -- they say --, let it be 
limited, but it should be perfectly grasped. Well, I am a mild and indulgent spirit, but 
every time this statement has been objected to me, a malicious thinking has crossed my mind as 
a flash: Oh, if only I could take my discussant at his word, and with a magic power I could 
turn off in his brain for one moment every vague cognition, by letting there only what he knows 
perfectly! Provided after such a cruel mutilation some glimmer would survive in his intellect, 
which I strongly doubt, this would resemble a will o' wisp lost into a deep and immense darkness. 
The truth is that we perfectly know nothing...}\footnote{G. Castelnuovo, 
{\it La scuola nei rapporti con la vita e la scienza moderna}, talk held in Genova in 1912 
at the III Congress of Mathesis, and reproduced in {\it Archimede}, n. 2--3, 1962.} 
\end{quote}

\vspace{0.1cm}

For further links between the approach to mathematics of Guido and Emma Castelnuovo we refer to \cite{F}.

\end{document}